\documentclass{article}
\usepackage{amsmath, amscd, amsfonts, amssymb, amsthm}
\newtheorem{theorem}{Theorem}[section]
\newtheorem{corollary}[theorem]{Corollary}
\newtheorem{lemma}[theorem]{Lemma}
\newtheorem{prop}[theorem]{Proposition}
\begin{document}

\title{Extremal Reversible Measures for the Exclusion Process\footnote{Research supported in part by NSF grant DMS-00-70465}}
\author{Paul Jung\footnote{Department of Mathematics, University of California, Los
Angeles, CA 90095-1555, USA}}

\maketitle

\begin{abstract}The invariant measures $\mathcal{I}$ for the exclusion
process have long been studied and a complete description is known
in many cases. This paper gives characterizations of $\mathcal{I}$
for exclusion processes on $\mathbb{Z}$ with certain reversible
transition kernels.  Some examples for which $\mathcal{I}$ is
given include all finite range kernels that are asymptotically
equal to $p(x,x+1)=p(x,x-1)=1/2$. One tool used in the proofs
gives a necessary and sufficient condition for reversible measures
to be extremal in the set of invariant measures, which is an
interesting result in its own right. One reason that this
extremality is interesting is that it provides information
concerning the domains of attraction for reversible measures.
\end{abstract}


\section{Introduction}
Given a countable set $\mathcal{S}$ and a corresponding
probability transition function $p(x,y)$ satisfying $\sup_y \sum_x
p(x,y)<\infty$, IPS (Liggett(1985)) constructs and describes the
exclusion process on $\{0,1\}^\mathcal{S}$. Its generator is given
by the closure of the operator $\Omega$ on
$\mathcal{D}(\{0,1\}^{\mathcal{S}})$, the set of all functions on
$\{0,1\}^\mathcal{S}$ that depend on finitely many coordinates. If
$f\in \mathcal{D}(\{0,1\}^{\mathcal{S}})$ and $\eta_{xy}$ is
defined as
\begin{equation*}
\eta_{xy}(u)=\left\{
\begin{array}{ll}
\eta(y)&\text{if }u=x\\
\eta(x)&\text{if }u=y\\
\eta(u)&\text{if }u\neq x,y\\
\end{array}
\right.
\end{equation*}
then
\begin{equation*}
{\Omega}f(\eta)=\sum_{\eta(x)=1,\eta(y)=0}{p(x,y)[f(\eta_{xy})-f(\eta)]}.
\end{equation*}
The semigroup of this process will be denoted by $S(t)$.

When $p(x,y)=p(y,x)$ the process has been completely studied in
that a full description of its invariant measures is known as well
as their respective domains of attraction. The asymmetric
exclusion process on the other hand has been much more elusive.
General classes of invariant measures are known in the two cases
where $p(x,y)$ is doubly stochastic (i.e., $\sum_{x\in\mathcal{S}}
p(x,y)=1$ for all $y\in\mathcal{S}$) or when there exists a
reversible measure $\pi(x)>0$ on $\mathcal{S}$ (i.e., a measure
satisfying $\pi(x)p(x,y)=\pi(y)p(y,x)$). However, a complete
description of $\mathcal{I}$ is known only in the three cases when
either
\begin{enumerate} \item[(a)] $p(x,y)$ is reversible and
positive recurrent for either the particles or holes ($1$'s or
$0$'s) (Liggett(1976)) \item[(b)] $p(x,y)$ corresponds to certain
random walks on $\mathbb{Z}$ (Liggett(1976) and Bramson, Liggett,
and Mountford(2002)) or \item[(c)] $p(x,y)$ corresponds to a birth
and death chain on $\mathbb{Z}^+$ (Liggett(1976)).
\end{enumerate}
Almost nothing is known about the domains of attraction concerning
invariant measures in the asymmetric case, although we note here
that there are some nice theorems concerning the case where
$p(x,y)$ is an asymmetric simple random walk on $\mathbb{Z}$ (see
Liggett(1999)).

Our purpose in this paper is to shed some more light on the
problem of classifying $\mathcal{I}$ and its respective domains of
attraction for the asymmetric exclusion process when a reversible
measure $\pi(x)$ exists for $p(x,y)$. In order to describe the
results of this paper we must first discuss case (a) and state a
special case of (b) above.

We start by stating what is known for the mean zero case of (b).
Let $\nu_\rho$ be the product measure on $\{0,1\}^{\mathcal{S}}$
with marginals $\nu_\rho\{\eta:\eta(x)=1\}=\rho$. Liggett(1976)
uses a coupling of two exclusion processes to show that when
$p(x,y)=p(0,y-x),\sum_x{|x|p(0,x)}<\infty, \text{ and }
\sum_x{xp(0,x)}=0$ on $\mathbb{Z}$ the set of extremal invariant
measures is
\begin{equation}\label{lthm}
\mathcal{I}_{e}=\{\nu_\rho:0\le\rho\le 1\}.
\end{equation}

Before describing the invariant measures for case (a), we define
some extremal reversible invariant measures $\{\nu^{(n)}\}$ when a
reversible measure $\pi(x)$ satisfying
\begin{equation}\label{<inf}
\sum_x{\pi(x)/[1+\pi(x)]^2}<\infty
\end{equation}
exists. This family of extremal reversible measures was first
discovered by Liggett.  In particular, he breaks down (\ref{<inf})
into three cases and writes
\begin{enumerate}
\item If $\sum_x \pi(x)<\infty$, let $A_n=\{\eta:\sum_x
\eta_x=n\}$ for nonnegative integers $n$.

\item If $\sum_x 1/\pi(x)<\infty$, let $A_n=\{\eta:\sum_x
[1-\eta_x]=n\}$ for nonnegative integers $n$.

\item If $\sum_x \pi(x)/[1+\pi(x)]^2<\infty$, $\sum_x
\pi(x)=\infty$, and $\sum_x 1/\pi(x)=\infty$, there exists a
$T\subset S$ for which $\sum_{x\in T} \pi(x)<\infty$ and
$\sum_{x\notin T} 1/\pi(x)<\infty$.  In this case, let
\begin{eqnarray*}
A_n=\{\eta:\sum_{x\in T} \eta(x)-\sum_{x\notin T}[1-\eta(x)]=n\}
\end{eqnarray*}
 for integers $n$.
\end{enumerate}

To define $\{\nu^{(n)}\}$, let $\nu^{c}$ be the product measure
with marginals
$\nu^c\{\eta:\eta(x)=1\}=\frac{c\pi(x)}{1+c\pi(x)}$. Liggett shows
that the measures
\begin{eqnarray*}\label{(n)}
\nu^{(n)}(\cdot)=\nu^c(\cdot|\mathcal{A}_n)\text{ for }n\in\mathbb{Z},\\
\nu^{(\infty)}=\text{ the pointmass on }\eta(x)\equiv 1,\nonumber\\
\nu^{(-\infty)}=\text{ the pointmass on }\eta(x)\equiv 0\nonumber
\end{eqnarray*}
are the unique stationary distributions for the positive recurrent
Markov chains on $A_n$. A simple consequence of Theorem B52 in
Liggett(1999) is that the reversible measures $\{\nu^{(n)}\}$ are
extremal in the set of invariant measures. For the first two cases
in the trichotomy of (\ref{<inf}) above, these are the only
extremal invariant measures. These first two cases correspond
exactly to (a) above. Note that changing $T$ in the third case of
the trichotomy above amounts to a relabeling of the sequence
$\{\nu^{(n)}, n\in\mathbb{Z}\}$.

Whenever a reversible measure $\pi(x)$ on $\mathcal{S}$ exists,
the product measures $\{\nu^{c}\}$ are well-defined. Theorem
VIII.2.1 in IPS tells us that these measures are invariant for the
exclusion process. Applying \mbox{Kakutani's} Dichotomy (e.g. page
244 of Durrett(1996))  we have that
$\sum_x{\pi(x)/[1+\pi(x)]^2}=\infty$ is a necessary and sufficient
condition for the measures \mbox{$\{\nu^{c}:0\le c\le\infty\}$} to
be mutually singular.
 Since all the results in this paper concern
the reversible measures $\{\nu^c\}$, we will assume throughout the
rest of the paper that $\pi(x)$ satisfying
$\pi(x)p(x,y)=\pi(y)p(y,x)$ exists.

In Section \ref{2} we prove Theorem \ref{lemma1.2} which states
that $\sum_x{\pi(x)/[1+\pi(x)]^2}=\infty$ is exactly the situation
in which the measures $\nu^c$ are extremal invariant. Not only
does this result have some nice applications, but knowing that an
invariant measure is extremal in the set of invariant measures has
always been an interesting question concerning particle systems.
Examples of such results are Theorem III.1.17 in Liggett(1999) and
Theorem 1.4 in Sethuraman(2001). The main reason extremality of
invariant measures is interesting is its close connection with
ergodicity.  This is seen by the application Theorem III.1.17 in
Liggett(1999) to prove Theorem III.4.8 in Liggett(1999) concerning
the tagged particle process; it is again seen by the application
of Theorem 1.4 in Sethuraman(2001) to certain central limit
theorems given in Kipnis and Varadhan(1986). In particular, if the
initial measure for a process is an extremal invariant measure
then the process evolution is ergodic with respect to time shifts.

Sections \ref{liggett} and \ref{z} use Theorem \ref{lemma1.2} to
extract information about the invariant measures of the process on
$\mathbb{Z}$.  In particular, Section \ref{liggett} modifies
Liggett's original proof of the result stated above equation
$(\ref{lthm})$ to obtain the following result:
\begin{theorem} \label{thm1.1}
Let $\mathbb{Z}$ be irreducible with respect to a transition
kernel $p(x,y)$ for which there exists a reversible measure
$\pi(x)$. Suppose $q_i(z)$ is a transition kernel such that
$\sum_{z}{zq_i(z)}=0$ and $\sum_{z}{|z|q_i(z)}<\infty$ for
$i=1,2$, and suppose that
\begin{equation}\label{maincond}
\lim_{K\rightarrow\infty}\sum_{x\ge 0}\sum_{|z|\ge|x-K|}
|p(x,x+z)-q_1(z)|=0 \text{ and }
\lim_{K\rightarrow\infty}\sum_{x\le 0}\sum_{|z|\ge|x+K|}
|p(x,x+z)-q_2(z)|=0.
\end{equation}

(a) If $\sum_x \pi(x)/[1+\pi(x)]^2=\infty$ then
$\mathcal{I}_e=\{\nu^{c}:0\le c\le\infty\}$.

(b) If $\sum_x \pi(x)/[1+\pi(x)]^2<\infty$ then
$\mathcal{I}_e=\{\nu^{(n)}\}$.
\end{theorem}

In essence the above theorem says that when the transition
probabilities are asymptotically translation invariant and have an
asymptotic mean of zero, the reversible measures are the only
invariant measures.  Theorem \ref{thm1.1} is merely an extension
(in the case where $\pi(x)$ exists) of the theorem proved by
Liggett(1976) which is stated above equation $(\ref{lthm})$.

Condition (\ref{maincond}) may seem somewhat daunting, but note
that if $\lim_{x\rightarrow\infty} p(x,x+z)=q_1(z)$,
$\lim_{x\rightarrow -\infty} p(x,x+z)=q_2(z)$, and $p(x,y)$ has
finite range (i.e. $p(x,y)\equiv 0$ if $|x-y|>n$ for some $n$),
then $(\ref{maincond})$ and $\sum_z|z|q_1(z)<\infty$ are both
automatically satisfied. Also, the below condition which is
somewhat easier to grasp than $(\ref{maincond})$ implies
$(\ref{maincond})$:
\begin{equation*}
\sum_{x\ge 0}\sum_{z} |p(x,x+z)-q_1(z)|<\infty \text{ and }
\sum_{x\le 0}\sum_{z} |p(x,x+z)-q_2(z)|<\infty.
\end{equation*}
A typical situation for which the theorem holds is when the
transition rates are nearest-neighbor and are given by
$p(x,x+1)=p(x,x-1)=1/2$ except for finitely many $x$.

Note that the premises of the theorem together with the assumption
that a reversible $\pi(x)$ exists imply that $q_i(z)$ must be
symmetric. To see this suppose $q_1(z)$ is not symmetric. Also,
assume that $q_1(z_1)>q_1(-z_1)>0$ for some $z_1\in\mathbb{N}$. We
can do this without loss of generality since $q_1(z)>0$ implies
$q_1(-z)>0$ by the reversibility of $\pi(x)$.  The mean zero
assumption tells us there exists $z_2\in\mathbb{N}$ such that
$q_1(z_2)<q_1(-z_2)$. If $z_3$ is a multiple of both $z_1$ and
$z_2$ then since $p(x,x+z)\rightarrow q_1(z)$ we can find $x_1$ so
that for $x>x_1$, $\pi(x)<\pi(x+z_3)$.  But we can also find $x_2$
so that for $x>x_2$, $\pi(x)>\pi(x+z_3)$, a contradiction. So
$q_1(z)$ must be symmetric. The proof that $q_2(z)$ is symmetric
follows similarly.

The proof of the above theorem follows Liggett's original outline
and does not actually require Theorem \ref{lemma1.2}.  However,
the usefulness of Theorem \ref{lemma1.2} is seen in the
simplification of one part of Liggett's original proof.

In Section \ref{z} we prove a theorem concerning the
nearest-neighbor exclusion process on $\mathbb{Z}$. For the
statement of the theorem we will need the following definitions.

Let $\mathcal{L}^-$ be the set of limit points of $\{\pi(x),x<0\}$
and $\mathcal{L}^+$ be the set of limit points of
$\{\pi(x),x>0\}$.
\begin{theorem}\label{thm1.8}
Suppose that $\inf_{|x-y|=1}p(x,y)>0$ for a nearest-neighbor
exclusion process on $\mathbb{Z}$. Then nonreversible invariant
measures can exist only when either (a) $\mathcal{L}^-=\{0\}$ and
$\mathcal{L}^+=\{\infty\}$ or (b) $\mathcal{L}^-=\{\infty\}$ and
$\mathcal{L}^+=\{0\}$.
\end{theorem}

The above theorem in no way guarantees the existence of
nonreversible invariant measures as seen by the following example.
Let
\begin{eqnarray}\label{example}
p(-1,-2)=p(-1,0)=p(0,-1)=p(0,1)=1/2,\\
p(x,x+1)=1-p(x,x-1)=\frac{|x|+1}{2|x|}\text{ otherwise}.\nonumber
\end{eqnarray}
This transition gives us situation (a) in the theorem above. The
reversible invariant measures $\{\nu^c\}$ certainly exist, but it
is easy to see that condition (b) of Theorem \ref{thm1.1} is
satisfied by (\ref{example}), therefore there are no nonreversible
invariant measures.

A curious aside is as follows. If in this example we start this
process off with initial measure $\nu_{\rho}$ and take the limit
of some converging sequence of measures
\begin{equation}\label{choose}
\frac{1}{T_n}\int_0^{T_n} \nu_{\rho}S(t)\, dt
\end{equation}
then Theorem I.1.8 in IPS says that this limit is an invariant
measure for the process. In view of the previous discussion, this
limit must converge to some mixture of the extremal invariant
measures $\{\nu^{(n)}, 0\le n\le\infty\}$. It would be interesting
indeed to find out which mixture (\ref{choose}) converges to. Note
here that we started off with an initial state that concentrates
on an uncountable number of states, but the limiting distribution
concentrates on a countable number of states (which may very well
be just the point masses of all $0$'s and all $1$'s).

If $p(x,y)$ is an asymmetric random walk kernel with nonzero mean
then we have one of the situations described in the theorem above,
and one might correctly guess that there exists some nonreversible
invariant measure. In fact, a well-known result of Liggett(1976)
proves that the measures
\begin{equation}\label{rhomeasures}
\{\nu_\rho:0\le\rho\le 1\}
\end{equation}
are invariant measures. Since any limit of (\ref{choose}) is
invariant, we intuitively might have expected this. More
precisely, if there were no nonreversible measures then this limit
would presumably be a mixture of the reversible measures $\nu^c$.
But it is intuitive that there is no mixture of $\nu^c$'s to which
this limit could converge, leading us to believe that the limit
converges to some other measure.

We note here that the set of measures in (\ref{rhomeasures}) is
the same as the set of measures in (\ref{lthm}) but are of an
entirely different nature.  In the setting of (\ref{lthm}) the
measures $\{\nu_\rho:0\le\rho\le 1\}$ are reversible and
constitute the entire set of extremal invariant measures.  On the
other hand, under the current setting, the measures
$\{\nu_\rho:0\le\rho\le 1\}$ are not reversible and
\begin{equation*}
\mathcal{I}_e= \{\nu_\rho:0\le\rho\le 1\}\cup\{\nu^c:0\le
c\le\infty\}.
\end{equation*}

The discussion in the previous paragraphs might make us wonder for
which transition kernels a nonreversible invariant measure exists.
To gain more insight into the situation we introduce a concept
known as the \textit{flux} of an invariant measure $\mu$. We will
continue to assume that the transition probabilities are
nearest-neighbor, but we will no longer assume they are
translation invariant. Define
\begin{equation}\label{flux}
\text{flux}(\mu)=p(x,x+1)\mu\{\eta:\eta(x)=1,
\eta(x+1)=0\}-p(x+1,x)\mu\{\eta:\eta(x)=0, \eta(x+1)=1\}.
\end{equation}
Let $1_x(\eta)=\eta(x)$ be the indicator function of
$\{\eta(x)=1\}$. By computing the positive and negative terms of
the left-hand side of $\int \Omega 1_x d\mu=0$ it can be seen that
$\text{flux}(\mu)$ is independent of $x$.

When an invariant measure $\mu$ is reversible it can easily be
seen from (\ref{flux}) that $\text{flux}(\mu)=0$. So if an
invariant measure exists whose flux is nonzero it must be
nonreversible. For the process with $p(x,x+1)>1/2$ and
$p(x,x-1)=1-p(x,x+1)$, the invariant measures
$\{\nu_\rho:0\le\rho\le 1\}$ all have a positive flux with the
flux being maximized when $\rho=1/2$ (a full discussion of this
can be found in either Janowski and Lebowitz(1994) or Part III of
Liggett(1999)).  This positive flux is the reason why
(\ref{rhomeasures}) is fundamentally different from (\ref{lthm}).
It would be quite nice if one could prove that some nonreversible
invariant measure exists whenever $p(x,x+1)>1/2+\epsilon$ for all
$x$. The $\epsilon$ here serves the role of providing some
positive flux in the limit.

Finally, Section \ref{domains} will apply Theorem \ref{lemma1.2}
to give information concerning the domains of attraction (in the
Cesaro sense) of reversible measures in the case where
$\sum_x{\pi(x)/[1+\pi(x)]^2}=\infty$. The results of Section
\ref{domains} only give sufficient conditions for Cesaro
convergence to an invariant measure, but are nonetheless
interesting since so little is known concerning domains of
attraction for the asymmetric exclusion process.  The key known
results concerning domains of attraction of asymmetric exclusion
processes are stated in Andjel, Bramson, Liggett(1988). They
concern the limiting distribution of exclusion processes with
asymmetric nearest-neighbor random walk kernels when the initial
measures are certain product measures. To get an idea of how
difficult it is to prove anything of this sort, we refer the
reader to Andjel, Bramson, Liggett(1988).

The fact that Theorem \ref{thmlastsection} concerns Cesaro
convergence rather than the usual weak convergence, while
undesirable, is not so bad since many results in particle systems
concern Cesaro convergence (see Section I.1 in IPS).  One notable
example of this is the main result of Andjel(1986) which concerns
the Cesaro convergence of certain initial product measures when
the transition kernel of the exclusion process is an asymmetric
nearest-neighbor random walk.  In fact these results were later
shown to be true for weak convergence (this was the goal of
Andjel, Bramson, Liggett(1988)).  We note here that Theorem
\ref{thmlastsection} does not use the property of reversibility,
therefore one can apply the theorem to situations in which one
knows that a particular invariant measure is extremal in the set
of invariant measures.

\section{Extremal reversible measures}\label{2}
In this section we state and prove Theorem \ref{lemma1.2}. The
common technique used in the proof of this theorem and in the
proofs of most of the other results in this paper is the coupling
technique. We now define the basic coupling of $\eta_t$ and
$\xi_t$ which lets the two exclusion processes move together as
much as possible. The generator for this coupling is the closure
of the operator $\tilde{\Omega}$ defined on
$\mathcal{D}(\{0,1\}^{\mathcal{S}}\times\{0,1\}^{\mathcal{S}})$:
\begin{eqnarray*}
&&\tilde{\Omega}f(\eta,\xi)=\sum_{\eta(x)=\xi(x)=1,\eta(y)=\xi(y)=0}{p(x,y)[f(\eta_{xy},\xi_{xy})-f(\eta,\xi)]}\\
&+&\sum_{\eta(x)=1,\eta(y)=0 \text{ and } (\xi(y)=1 \text{ or } \xi(x)=0)}{p(x,y)[f(\eta_{xy},\xi)-f(\eta,\xi)]}\\
&+&\sum_{\xi(x)=1,\xi(y)=0 \text{ and } (\eta(y)=1 \text{ or }
\eta(x)=0)}{p(x,y)[f(\eta,\xi_{xy})-f(\eta,\xi)]}.
\end{eqnarray*}

\begin{theorem}\label{lemma1.2}
Suppose $\mathcal{S}$ is irreducible with respect to $p(x,y)$.
Then the measures $\nu^c$ are extremal invariant if and only if
$\sum_x{\pi(x)/[1+\pi(x)]^2}=\infty$.
\end{theorem}

\begin{proof}
The discussion on page 383 of IPS shows that if
$\sum_x{\pi(x)/[1+\pi(x)]^2}<\infty$ then the measures $\nu^c$ are
not extremal invariant giving us one direction of the theorem.  We
will prove the other direction.

Assume throughout that $0<c<\infty$.  Since the measures $\nu^c$
are invariant and since all bounded continuous functions can be
approximated uniformly by functions that depend on finitely many
coordinates then by Theorem B52 in Liggett(1999), we need only
show that for any two functions $f$ and $g$ which depend on
finitely many coordinates
\begin{equation*}
\lim_{T\rightarrow\infty}\frac{1}{T}\int_{0}^{T}E^{\nu^c}f(\eta_0)g(\eta_t)dt=\int
f \,d\nu^c\int g\, d\nu^c.
\end{equation*}

We claim that to show the above equation holds, it is enough to
show that for any finite $A\subset\mathcal{S}$ and for
$\mu_{1,A}^c(\cdot)=\nu^c(\cdot|\{\eta:\eta(x)=1 \forall x\in
A\})$
\begin{equation}\label{equ3}
\lim_{T\rightarrow\infty}\frac{1}{T}\int_0^T \mu_{1,A}^cS(t)
dt=\nu^c.
\end{equation}
To see this define the measures
$\mu_{\zeta,A}^c(\cdot)=\nu^c(\cdot|\{\eta(x)=\zeta(x) \forall
x\in A\})$ where $\zeta$ is a configuration on $\{0,1\}^A$. We can
write the measure $\nu^c$ as a linear combination
\begin{equation*}
\nu^c=\sum_{\zeta\in \{0,1\}^A}{a_{\zeta}\mu_{\zeta,A}^c}
\end{equation*}
where we use the convention that $\zeta=i$ is the configuration in
$\{0,1\}^A$ such that $\zeta(x)=i$ for all $x\in A$. For
\begin{equation*}
f_{A}=\left\{
\begin{array}{ll}
1   &\text{when }\eta(x)=1 \text{ for all }x\text{ in
 the finite set }A\\
0   &\text{otherwise}
\end{array}
\right.
\end{equation*}
we have that
\begin{equation*}
\lim_{T\rightarrow\infty}\frac{1}{T}\int_0^T
E^{\nu^c}f_A(\eta_0)g(\eta_t)dt
=\lim_{T\rightarrow\infty}\frac{1}{T}\int_0^Ta_{1}\int
S(t)g(\eta)d\mu_{1,A}^c dt=\int f_A d\nu^c\int g\, d\nu^c
\end{equation*}
which proves the claim.

Define $\mu_{0,A}^c$ similarly to the way we defined
$\mu_{1,A}^c$.  If we assume a fixed $A$ then we can drop the
subscript $A$ so as to write $\mu_i^c=\mu_{i,A}^c$. The rest of
the proof will now argue that (\ref{equ3}) holds.

Choose $\delta>0$ and couple the processes $\eta_t$ and $\xi_t$
using the basic coupling starting with measures $\mu_0^c$ and
$\mu_1^{c}$ so that $\eta_0$ and $\xi_0$ disagree only for $x\in
A$. In particular, since the basic coupling is the coupling which
allows $\eta_t$ and $\xi_t$ to move together as much as possible,
then $\eta_t$ and $\xi_t$ can differ at most at $n$ sites where
$|A|=n$.

If there exists $\bar{T}$ such that for all $T>\bar{T}$
\begin{equation*}
\frac{1}{T}\int_0^{T}[\mu_1^cS(t)\{\xi(0)=1\}-
\mu_0^cS(t)\{\eta(0)=1\}]dt\le \delta
\end{equation*}
then we must have that
\begin{equation*}
\lim_{T\rightarrow\infty}\frac{1}{T}\int_0^{T}\mu_1^cS(t)\{\xi(0)=1\}dt=
\lim_{T\rightarrow\infty}\frac{1}{T}\int_0^{T}\mu_0^cS(t)\{\eta(0)=1\}dt.
\end{equation*}
Keeping in mind the way that $\eta_t$ and $\xi_t$ are coupled,
irreducibility then tells us that
\begin{equation*}
\lim_{T\rightarrow\infty}\frac{1}{T}\int_0^{T}\mu_1^cS(t)dt=
\lim_{T\rightarrow\infty}\frac{1}{T}\int_0^{T}\mu_0^cS(t)dt.
\end{equation*}
But the measure $\nu^c$ lies stochastically between the left-hand
side and the right-hand side of the equation above, so in fact we
must have that (\ref{equ3}) holds.

We can therefore assume to the contrary that there exists a
$\delta>0$ and a sequence $\{T_n\}$ such that
\begin{equation}\label{eqndelta}
\frac{1}{T_n}\int_0^{T_n}[\mu_1^cS(t)\{\xi(0)=1\}-
\mu_0^cS(t)\{\eta(0)=1\}]dt>\delta
\end{equation}
for all $n$.

Pick $\epsilon>0$ so that
\begin{equation*}
\nu^{c+\epsilon}\{\xi(0)=1\}-\nu^{c-\epsilon}\{\eta(0)=1\}<\delta/3.
\end{equation*}
Using the basic coupling once more, couple the processes $\eta_t$
and $\xi_t$ starting off in the measures $\mu_1^c$ and
$\nu^{c+\epsilon}$ so that $\lambda_1\{(\eta,\xi):\eta(x)\le
\xi(x) \text{ for all }x\in \mathbb{Z}\backslash A\}=1$ where
$\lambda_1$ is the coupling measure. If
$\hat{\mu}^c=\nu^c(\cdot|\{\eta:\eta(x)=0 \text{ for some } x\in
A\})$ then
\begin{equation*}
\nu^c=\gamma\mu_1^c+(1-\gamma)\hat{\mu}^c
\end{equation*}
for $\gamma=\nu^c\{\eta:\eta(x)=1 \forall x\in A\}$. Couple the
measures $\hat{\mu}^c$ and $\nu^{c+\epsilon}$ in a way similar to
$\lambda_1$ so that we get another coupling measure $\lambda_2$.

Choose a subsequence $\{T_{n_k}\}$ so that we can define some
limiting invariant measure
\begin{equation*}
\omega_1=\lim_{k\rightarrow\infty}\frac{1}{T_{n_k}}\int_0^{T_{n_k}}\lambda_1
S(t)dt.
\end{equation*}
Let $\nu_1^c$ be the $\eta$-marginal measure of $\omega_1$ so that
in particular
\begin{equation*}
\nu_1^c=\lim_{k\rightarrow\infty}\frac{1}{T_{n_k}}\int_0^{T_{n_k}}\mu_1^cS(t)dt.
\end{equation*}

To complete the proof of the theorem we will need the following
lemma:
\begin{lemma}\label{interlemma}
$\nu^{c+\epsilon}\ge\nu_1^c$.
\end{lemma}

\begin{proof}[Proof of lemma]
Let $f_x(\eta,\xi)=[1-\eta(x)]\xi(x)$, $ D_m=\{(\eta,\xi):\eta(x)>
\xi(x)\text{ at exactly }m\text{ sites}\}$, and $D=\bigcup_{m\ge
1} D_m$. If $\nu^{c+\epsilon}\ngeq\nu_1^c$ then it must be that
$\omega_1(D)>0$. We claim that this implies
\begin{equation*}
\int_{D}\sum_{x} f_x d\omega_1=0.
\end{equation*}

To prove the claim, assume to the contrary that $\int_{D}\sum_{x}
f_x d\omega_1>0$ so that there exist sites for which
$\eta(x)<\xi(x)$. Let $M$ be the largest $m$ for which
$\omega_1(D_m)>0$. Then by the irreducibility condition and by the
fact that there exist sites for which $\eta(x)<\xi(x)$ we have
$\omega_1S(t)(D_M)<\omega_1(D_M)$ for $t>0$. But this is a
contradiction to the invariance of $\omega_1$ proving the claim.

Now if the two processes $\eta_t$ and $\xi_t$ have the measures
$\nu^c$ and $\nu^{c+\epsilon}$ respectively then let $\omega$ be
the coupling measure for $\{(\eta_t,\xi_t)\}$ which concentrates
on $\nu^c\le\nu^{c+\epsilon}$. For this coupling, the $\omega$
probability that $f_x(\eta,\xi)=1$ for a given $x$ is equal to the
left-hand side below:
\begin{equation*}
\frac{(c+\epsilon)\pi(x)}{1+(c+\epsilon)\pi(x)}-\frac{c\pi(x)}{1+c\pi(x)}>
\frac{\epsilon\pi(x)}{[1+(c+\epsilon)\pi(x)]^2}.
\end{equation*}
Since $\sum_x{\pi(x)/[1+\pi(x)]^2}=\infty$, by the Borel-Cantelli
Lemma the $\omega$ probability that $\sum_{x} f_x=\infty$ is equal
to $1$. The measure $\omega_1$ is absolutely continuous with
respect to $\omega$ since
\begin{equation*}
\omega=\gamma\lambda_1 +
(1-\gamma)\lambda_2=\gamma\omega_1+(1-\gamma)\lim_{k\rightarrow\infty}\frac{1}{T_{n_k}}\int_0^{T_{n_k}}
\lambda_2S(t)dt
\end{equation*}
where $\lambda_2$ is as defined above. Therefore $\int_{E}\sum_{x}
f_x d\omega_1=\infty$ for any set $E$ with positive $\omega_1$
measure which contradicts $\int_{D}\sum_{x} f_x d\omega_1=0$ so it
must be that $\omega_1(D)=0$ proving the lemma.
\end{proof}

We now turn back to the proof of the theorem. Since by the lemma
we have $\nu^{c+\epsilon}\ge\nu_1^c$, then there exists a $K$ such
that for all $k>K$
\begin{equation*}
\frac{1}{T_{n_k}}\int_0^{T_{n_k}}\mu_1^cS(t)\{\eta(0)=1\}dt-\nu^{c+\epsilon}\{\xi(0)=1\}
<\delta/3.
\end{equation*}
If $\nu_0^c$ is some limiting measure of
\begin{equation*}
\frac{1}{T_{n_{k_l}}}\int_0^{T_{n_{k_l}}}\mu_0^cS(t)dt
\end{equation*}
then an argument similar to that used in Lemma \ref{interlemma}
shows that $\nu^{c-\epsilon}\le\nu_0^c$. There then exists an $L$
such that for $l>L$
\begin{equation*}
\nu^{c-\epsilon}\{\eta(0)=1\}-\frac{1}{T_{n_{k_l}}}\int_0^{T_{n_{k_l}}}
\mu_0^{c}S(t)\{\xi(0)=1\}dt<\delta/3.
\end{equation*}
Altogether we have for $l>L$,
\begin{equation*}
\frac{1}{T_{n_{k_l}}}\int_0^{T_{n_{k_l}}}[\mu_1^cS(t)\{\xi(0)=1\}-\mu_0^{c}S(t)\{\eta(0)=1\}]dt<\delta
\end{equation*}
which contradicts inequality (\ref{eqndelta}) so it must be that
(\ref{equ3}) holds completing the proof of the theorem.
\end{proof}

\section{The asymptotically mean zero process on $\mathbb{Z}$}\label{liggett}

In this section we prove Theorem \ref{thm1.1}.  To do so we will
need to define $\mathcal{\tilde{I}}$ as the set of invariant
measures for the basic coupling and $\mathcal{\tilde{I}}_{e}$ as
its extreme points.

Recall that $
 f_x(\eta,\xi)=[1-\eta(x)]\xi(x)$. In order to simplify the notation we further define the
functions
\begin{equation*}
\begin{array}{ll}h_{yx}(\eta,\xi)=[1-\eta(y)][1-\xi(y)]f_x(\eta,\xi),
&g_{yx}(\eta,\xi)=\eta(y)\xi(y) f_x(\eta,\xi),\\\text{and
}f_{yx}(\eta,\xi)=\eta(y)[1-\xi(y)]f_x(\eta,\xi).
\end{array}
\end{equation*}
In particular, for $T$ a finite subset of $\mathcal{S}$ we have
\begin{eqnarray}\label{1}
&&\tilde{\Omega}\left(\sum_{x\in T} {f_x(\eta,\xi)}\right)
=-\sum_{x\in {T},y\in\mathcal{S}}{(p(x,y)+p(y,x))f_{yx}(\eta,\xi)}\\
&+&\sum_{x \in {T}, y \notin
{T}}{\left[p(x,y)g_{xy}-p(y,x)g_{yx}\right]} +\sum_{x\in {T},
y\notin{ {T}}}{\left[p(y,x)h_{xy}-p(x,y)h_{yx}\right]}\nonumber.
\end{eqnarray}

\begin{proof}[Proof of Theorem \ref{thm1.1}]
Let $\nu\in\mathcal{\tilde{I}}$. Then
$\int{\tilde{\Omega}(\sum_{x\in T}{f_x})}d\nu=0$ for each finite
${T}\subset \mathbb{Z}$ so that for
${T}_{[m,n]}=\{x\in\mathbb{Z}:m\le x\le{n}\}$ we get
\begin{eqnarray}\label{suc}
&&\sum_{x\in{T}_{[m,n]}, y\in\mathbb{Z}}{(p(x,y)+p(y,x))\int f_{yx} d\nu}\\
&=&\sum_{x \in {T}_{[m,n]}, y \notin {T}_{[m,n]}}{p(x,y)\int
(g_{xy} - h_{yx})d\nu} +\sum_{x\in{T}_{[m,n]},
y\notin{{T}}_{[m,n]}}{p(y,x)\int (h_{xy} - g_{yx}) d\nu}\nonumber.
\end{eqnarray}
 Notice that the left-hand side of this equation is increasing in
$n$ and $-m$, so that when we take the limit as
$n\rightarrow\infty$ or as $-m\rightarrow\infty$, a limit exists.

Choosing $\epsilon>0$ we can find $N$ so that for $n>N$:
\begin{eqnarray*}
&& \sum_{x>n+N}\sum_{z<n-x} p(x,x+z)\le\sum_{x>n+N}\sum_{z<n-x}
|p(x,x+z)-q_1(z)|+\sum_{|z|>N}{|z|q_1(z)}<\frac{\epsilon}{3}\\
&& \sum_{0<x<n}\sum_{z>n-x+N}
p(x,x+z)\le\sum_{0<x<n}\sum_{z>n-x+N}
|p(x,x+z)-q_1(z)|+\sum_{|z|>N}{|z|q_1(z)}<\frac{\epsilon}{3}\\
&&\sum_{ x\le 0}\sum_{z>n+N-x} p(x,x+z)\le\sum_{x\le
0}\sum_{z>n+N-x}
|p(x,x+z)-q_2(z)|+\sum_{|z|>N}{|z|q_2(z)}<\frac{\epsilon}{3}\\
\text{and}\\
&&\sum_{x<-n-N}\sum_{z>-x-n} p(x,x+z)\le\sum_{x<-n-N}\sum_{z>-x-n}
|p(x,x+z)-q_2(z)|+\sum_{|z|>N}{|z|q_2(z)}<\frac{\epsilon}{3}\\
&&\sum_{-n<x<0}\sum_{z<-x-n-N}p(x,x+z)\le\sum_{-n<x<0}\sum_{z<-x-n-N}
|p(x,x+z)-q_2(z)|+\sum_{|z|>N}{|z|q_2(z)}<\frac{\epsilon}{3}\\
&&\sum_{x\ge 0 }\sum_{z<-n-N-x}p(x,x+z)\le\sum_{x\ge
0}\sum_{z<-n-N-x}
|p(x,x+z)-q_1(z)|+\sum_{|z|>N}{|z|q_1(z)}<\frac{\epsilon}{3}.
\end{eqnarray*}
Since the construction of the exclusion process assumes that
$\sup_y \sum_x p(x,y)$ is finite (See IPS Chapter VIII) and since
$\int (g_{xy} - h_{yx})d\nu\le 1$, the right-hand side sums in
(\ref{suc}) above are absolutely convergent for any fixed $n$ and
$m$.

Now by the inequalities above and by $(\ref{maincond})$ we can
pass to the limit in (\ref{suc}) so as to write
\begin{eqnarray*}
&&\lim_{m\rightarrow
-\infty}\lim_{n\rightarrow\infty}\sum_{x\in{T}_{[m,n]},
y\in\mathbb{Z}}{(p(x,y)+p(y,x))\int f_{yx}
d\nu}\\
&=&\lim_{n\rightarrow\infty}\sum_{x \in {T}_{[0,n]},
y>n}{\left[q_1(y-x)\int (g_{xy}- h_{yx}) d\nu +
q_1(x-y)\int (h_{xy}- g_{yx}) d\nu\right]}\\
&+&\lim_{m\rightarrow -\infty}\sum_{x \in {T}_{[m,0]},
y<m}{\left[q_2(y-x)\int (g_{xy}- h_{yx}) d\nu + q_2(x-y)\int
(h_{xy}- g_{yx}) d\nu\right]}.
\end{eqnarray*}

The right-hand side above is equal to
\begin{eqnarray}\label{eq1.2}
&&\lim_{l\rightarrow\infty}\frac{1}{l}\sum_{n=1}^{l}\sum_{x \in
{T}_{[0,n]}, y>n}{\left[q_1(y-x)\int (g_{xy}- h_{yx}) d\nu +
q_1(x-y)\int (h_{xy}- g_{yx})
d\nu\right]}\\
&+&\lim_{k\rightarrow\infty}\frac{1}{k}\sum_{m=-1}^{-k}\sum_{x \in
{T}_{[m,0]}, y<m}{\left[q_2(y-x)\int (g_{xy}- h_{yx}) d\nu +
q_2(x-y)\int (h_{xy}- g_{yx}) d\nu\right].}\nonumber
\end{eqnarray}
We will devote the next few paragraphs to showing that these
limits are in fact equal to zero.

Define the measures $\nu^+$ and $\nu^-$ by choosing a subsequence
$n_j$ so that the following limits exist:
\begin{equation*}
  \nu^{+}=\lim_{j\rightarrow \infty}
  \frac{1}{n_j}\sum_{1\le{x}\le{n_j}}{\nu_{x}}
\end{equation*}
\begin{equation*}
  \nu^{-}=\lim_{j\rightarrow \infty}
  \frac{1}{|n_{-j}|}\sum_{-1\ge{x}\ge{n_{-j}}}{\nu_{x}}
\end{equation*}
where $\nu_x$ is the $x$ translate of $\nu$. In the partial sums
of (\ref{eq1.2}) above, for $j$ large enough each term
\begin{equation*}
q_i(y-x)\int(g_{xy}-h_{yx})d\nu
\end{equation*}
appears $|y-x|$ times when $q_i(y-x)>0$, so we can write
(\ref{eq1.2}) as
\begin{eqnarray} \label{eq1.3}
&&\sum_{z\in{\mathbb{Z}^+}}{[z q_1(z)\int (g_{oz}-h_{zo}) d\nu^+ -
z q_1(-z)\int (g_{zo}-h_{oz}) d\nu^+]}\\
&+&\sum_{z\in{\mathbb{Z}^-}}{[-z q_2(z)\int (g_{oz}-h_{zo}) d\nu^-
+ z q_2(-z)\int (g_{zo}-h_{oz}) d\nu^-]}\nonumber
\end{eqnarray}

Now consider two coupled processes with transition rates equal to
$q_1(z)$ and $q_2(z)$ respectively. The measures $\nu^+$ and
$\nu^-$ are translation invariant and are also invariant measures
for the coupled process with respect to $q_1(z)$ and $q_2(z)$
respectively. In particular if $\tilde{\Omega}_i$ is the generator
for the coupled process of $q_i(z)$, $\tilde{\Omega}$ is the
generator for the coupled process of $p(x,y)$, and
\begin{equation*}
f_{(A,B)}=\left\{
\begin{array}{ll}
1   &\text{when }\eta(x)=\xi(y)=1 \text{ for all }x\text{ in
 the finite set }A, y\text{ in the finite set }B\\
0   &\text{otherwise}
\end{array}
\right.
\end{equation*}
then
\begin{equation*}
\int\tilde{\Omega}_1 f_{(A,B)}
d\nu_+=\lim_{k\rightarrow\infty}\frac{1}{n_k}\sum_{1\le x\le
n_k}\int \tilde{\Omega}_1 f_{(A+x,B+x)}
d\nu=\lim_{k\rightarrow\infty}\frac{1}{n_k}\sum_{1\le x\le
n_k}\int \tilde{\Omega} f_{(A+x,B+x)} d\nu=0
\end{equation*}
where $A+x$ is the $x$ translate $A$.

By Lemma VIII.3.2 in IPS we have $\int f_{xy} d\nu^+=0$ for all
$x,y$. We can therefore write $\nu^+$ as $\nu^+=\lambda\nu_1
+(1-\lambda)\nu_2$ where $\nu_1$ concentrates on
$\{(\eta,\xi):\eta < {\xi}\}$ and $\nu_2$ on
$\{(\eta,\xi):\eta\ge{\xi}\}$.  Then
\begin{eqnarray*}
&&\int (g_{oz}-h_{zo}) d\nu^+=\lambda \int (g_{oz}-h_{zo})
d\nu_1\\
&=&\lambda[\nu_1\{(\eta,\xi):\eta(0)=1,
\eta(z)=0\}-\nu_1\{(\eta,\xi):\eta(0)=\xi(0)=1,\eta(z)=\xi(z)=0\}\\
&+&\nu_1\{(\eta,\xi):\eta(0)=\xi(0)=1,\eta(z)=\xi(z)=0\}-
\nu_1\{(\eta,\xi):\xi(0)=1,\xi(z)=0\}]\\
 &=&\lambda[\nu_1\{(\eta,\xi):\eta(0)=1,
\eta(z)=0\}-\nu_1\{(\eta,\xi):\xi(0)=1,\xi(z)=0\}].\\
\end{eqnarray*}

Because $\nu^+$ is translation invariant and invariant for the
process with rates $q_1(z)$, $\nu_1$ is also since $\nu_1$ and
$\nu_2$ are mutually singular and
$\nu^+=\lambda\nu_1+(1-\lambda)\nu_2$. By Theorem VIII.3.9 in IPS,
the marginals of $\nu_1$ are exchangeable, thus the right-hand
side above is equal to a constant $c^+$ as is the expression $\int
(g_{zo}-h_{oz}) d\nu^+$. Similarly we have that $\int
(g_{oz}-h_{zo}) d\nu^-$ and $\int (g_{zo}-h_{oz}) d\nu^-$ are
equal to a constant $c^-$.  Now by the mean zero assumption, we
have that expression (\ref{eq1.3}) is equal to $0$, but since
(\ref{eq1.2}) and (\ref{eq1.3}) are equal, we have in fact that
\begin{equation} \label{eq1.4}
\sum_{y\in{T}}{(p(x,y)+p(y,x))\int f_{xy} d\nu}=0
\end{equation}
for each finite $T\subset\mathbb{Z}$.

By irreducibility, if $\int f_{xy} d\nu>0$ for some $x,y$ then
$\int f_{xy} d\nu>0$ for all $x,y$. Choose $x_0$ and $y_0$ such
that $p(x_0,y_0)+p(y_0,x_0)>0$.  By (\ref{eq1.4}) and the
nonnegativity of $\int f_{xy} d\nu$, we must have that $\int
f_{x_0y_0} d\nu =0$ and thus $\int f_{xy} d\nu =0$ for all $x,y$.
Therefore $\nu\in \tilde{\mathcal{I}}$ implies that
\begin{equation*}
\nu\{(\eta, \xi):\eta< \xi \text{ or } \eta\ge \xi\}=1.
\end{equation*}

If $\sum_x{\pi(x)/[1+\pi(x)]^2}=\infty$ we can use Theorem
\ref{lemma1.2} to pick $\mu\in \mathcal{I}_e$ and $\nu^{c}\in
\mathcal{I}_e$. On the other hand if
$\sum_x{\pi(x)/[1+\pi(x)]^2}<\infty$ we can use the analysis in
the introduction to pick $\mu\in \mathcal{I}_e$ and $\nu^{(n)}\in
\mathcal{I}_e$. Since $\nu\{(\eta, \xi):\eta< \xi \text{ or }
\eta\ge \xi\}=1$, Proposition VIII.2.13 in IPS tells us there
exists a coupling with invariant measure $\nu$ where $\nu$ has
marginals $\mu\le\nu^c$ or $\mu\ge\nu^c$ in the first case and
marginals $\mu\le\nu^{(n)}$ or $\mu\ge\nu^{(n)}$ in the second
case.

Take first the case where $\sum_x{\pi(x)/[1+\pi(x)]^2}=\infty$.
Supposing that $\mu\neq \nu^0\neq\nu^\infty$, we have that there
exists a $c_0$ for which $\nu^{c_1}\le \mu$ for all $c_1<c_0$ and
$\mu\le\nu^{c_2}$ for all $c_2>c_0$. By the continuity of the one
parameter family of measures $\{\nu^c\}$ it must be that
$\mu=\nu^{c_0}$.

If $\sum_x{\pi(x)/[1+\pi(x)]^2}<\infty$ then we have three cases
(i), (ii), and (iii) as given in the introduction. Theorem
VIII.2.17 in IPS proves the first two cases so we will consider
only (iii). If $\mu\neq \nu^{(-\infty)}\neq\nu^{(\infty)}$ then
there exists an $n\in\mathbb{Z}$ such that either $\mu=\nu^{(n)}$
or $\nu^{(n)}<\mu<\nu^{(n+1)}$.  If the latter is true then $\mu$
concentrates on $A=\{\eta:\sum_{x\in T} \eta(x)<\infty,
\sum_{x\notin T} [1-\eta(x)]<\infty\}$ for some $T\subset S$ which
means that it must be a mixture of stationary distributions  for
the Markov chains on $A_n$ as described in the introduction.  But
$\mu\in\mathcal{I}_e$ so it must in fact be equal to some
$\nu^{(n)}$ completing the proof.
\end{proof}

We include in this section two more results which have proofs
similar to that of Theorem \ref{thm1.1}. We first need the
following definition: given transition probabilities $p(x,y)$
define the boundary of a set $\mathcal{T}$ to be
\begin{equation*}
\partial{\mathcal{T}}=\{x\notin\mathcal{T}:p(x,y)>0 \text{ for some }
y\in\mathcal{T}\}.
\end{equation*}

\begin{prop} \label{thm1.2}
Let $\mathcal{S}$ be irreducible with respect to $p(x,y)$ and
suppose that $\sum_x{\pi(x)/[1+\pi(x)]^2}=\infty$. If there exists
a sequence of increasing sets $\mathcal{T}_n$ such that
$\cup{\mathcal{T}_n}=\mathcal{S}$ and either
$\lim_{n\rightarrow\infty}\sum_{x\in{\partial\mathcal{T}_n}}{\pi(x)}=0$
or
$\lim_{n\rightarrow\infty}\sum_{x\in{\partial\mathcal{T}_n}}{1/\pi(x)}=0$,
then $\mathcal{I}_e=\{\nu^{c}:0\le c\le\infty\}$.
\end{prop}
\begin{proof}
Choose $\mu\in\mathcal{I}_e$. If
$\lim_{n\rightarrow\infty}\sum_{x\in{\partial\mathcal{T}_n}}{\pi(x)}=0$
then couple $\eta_t$ with $\xi_t$ so that they have the measures
$\mu$ and $\nu^c$ respectively.  If
$\lim_{n\rightarrow\infty}\sum_{x\in{\partial\mathcal{T}_n}}{1/\pi(x)}=0$
then couple them vice versa. We will prove the case in which
$\lim_{n\rightarrow\infty}\sum_{x\in{\partial\mathcal{T}_n}}{\pi(x)}=0$.
The other case follows similarly.

By (\ref{1}),
\begin{eqnarray*}
&&\sum_{x\in\mathcal{T}_{n}, y\in\mathcal{S}}{[p(x,y)+p(y,x)]\int f_{yx} d\nu}\\
&=&\sum_{x \in \mathcal{T}_{n}, y \notin \mathcal{T}_n}{p(x,y)\int
(g_{xy} - h_{yx})d\nu}
+\sum_{x\in\mathcal{T}_{n}, y\notin{\mathcal{T}}_n}{p(y,x)\int (h_{xy} - g_{yx}) d\nu}.\\
\end{eqnarray*}
Just as in the above proof, the left-hand side of this equation is
increasing in $n$ so that a limit exists as $n\rightarrow\infty$.
The right-hand side above goes to $0$ as $n\rightarrow\infty$
since
\begin{eqnarray*}
&&\sum_{x \in \mathcal{T}_{n}, y \notin \mathcal{T}_n}{p(x,y)\int
(g_{xy} - h_{yx})d\nu}
+\sum_{x\in\mathcal{T}_{n}, y\notin{\mathcal{T}}_n}{p(y,x)\int (h_{xy} - g_{yx}) d\nu}\\
&\le&\sum_{x \in \mathcal{T}_{n}, y \notin \mathcal{T}_n}{p(x,y)\int
f_yd\nu}
+\sum_{x\in\mathcal{T}_{n}, y\notin{\mathcal{T}}_n}{p(y,x)\int f_yd\nu}\\
&\le& C\sum_{y\in \partial\mathcal{T}_{n}}{\int
f_yd\nu}+\sum_{y\in \partial{T}_{n}}{\int f_yd\nu}\le
C\sum_{y\in
\partial\mathcal{T}_{n}}{\pi(y)}+\sum_{y\in
\partial{T}_{n}}{\pi(y)}.
\end{eqnarray*}
Here $C=\sup_y \sum_x p(x,y)$ which is finite by the assumptions
in the introduction.

Irreducibility now gives us $\int f_{xy} d\nu=0$ for all $x,y$.
The rest of the proof just follows that of Theorem \ref{thm1.1}.
\end{proof}

Note that if we change the hypothesis
$\sum_x{\pi(x)/[1+\pi(x)]^2}=\infty$ to
$\sum_x{\pi(x)/[1+\pi(x)]^2}<\infty$ then Theorem VIII.2.17 in IPS
says that $\mathcal{I}_e=\{\nu^{(n)}:0\le n\le\infty\}$.

\begin{corollary}\label{cor1.1}
If in Theorem \ref{thm1.1} we replaced condition
$(\ref{maincond})$ with the condition that $p(x,y)$ has finite
range, $\lim_{x\rightarrow+\infty}{p(x,x+z)}=q_1(z)$, and
$\lim_{x\rightarrow-\infty}\pi(x)$ equals $0$ or $\infty$ (or
alternatively $\lim_{x\rightarrow-\infty}{p(x,x+z)}=q_2(z)$, and
$\lim_{x\rightarrow+\infty}\pi(x)$ equals $0$ or $\infty$) then
the result still holds.
\end{corollary}

\begin{proof}
Replace expression (\ref{eq1.2}) in the proof of Theorem
\ref{thm1.1} with
\begin{eqnarray*}
&&\lim_{k\rightarrow\infty}\frac{1}{k}\sum_{n=1}^{k}\sum_{x \in
{T}_{[0,n]}, y>n}{\left[q_1(y-x)\int (g_{xy}- h_{yx}) d\nu +
q_1(x-y)\int (h_{xy}- g_{yx})
d\nu\right]}\\
&+&\lim_{m\rightarrow -\infty}\sum_{x \in {T}_{[m,0]},
y<m}\left[{p(x,y)\int (g_{xy} - h_{yx})d\nu} +{p(y,x)\int (h_{xy}
- g_{yx})d\nu}\right].
\end{eqnarray*} The proofs of Theorem \ref{thm1.1} and Proposition \ref{thm1.2} imply that this
expression is $0$.  The rest is proven above.
\end{proof}

Before moving on to the next section let us discuss what the above
results tell us in the case where $p(x,y)$ has finite range on
$\mathbb{Z}$. Proposition \ref{thm1.2} together with Theorem
VIII.2.17 in IPS says that if $\lim_{|x|\rightarrow\infty}\pi(x)$
equals $0$ or $\infty$ then the reversible measures are the only
invariant measures.  If the limits
$\lim_{x\rightarrow\infty}\pi(x)$ and
$\lim_{x\rightarrow-\infty}\pi(x)$ exist and one of them is
nonzero and finite, then the combination of Theorem \ref{thm1.1}
and Corollary \ref{cor1.1} imply that the only invariant measures
are the reversible ones. All together we have the following: if
$\pi(x)$ exists and has limits in both directions for the finite
range exclusion process on $\mathbb{Z}$, then unless the limit is
$0$ in one direction and $\infty$ in the other direction, the only
invariant measures are the reversible ones. Of course, as seen in
an example in the introduction, it is also possible to have
$\lim_{x\rightarrow\infty} p(x,x+z)=q_1(z)$ and
$\lim_{x\rightarrow -\infty} p(x,x+z)=q_2(z)$ as given in Theorem
\ref{thm1.1} and at the same time have the limit of $\pi(x)$ to be
$0$ in one direction, $\infty$ in the other.  In those cases
Theorem \ref{thm1.1} rules out nonreversible invariant measures. A
similar comment can be made for Corollary \ref{cor1.1}. We remind
the reader, however, that if the transition probabilities are
translation invariant with a drift so that the limit of $\pi(x)$
is $0$ in one direction and $\infty$ in the other direction, then
Liggett(1976) tells us that $\{\nu_\rho:0\le\rho\le1\}$ is a class
of nonreversible invariant measures.

\section{The nearest-neighbor process on $\mathbb{Z}$}\label{z}

We now restrict our attention to the nearest-neighbor case. More
specifically, assume throughout this section that we are dealing
with the irreducible nearest-neighbor exclusion process on
$\mathbb{Z}$ ($p(x,y)=0$ if and only if $|x-y|>1$).  In this case,
a reversible $\pi(x)$ always exists so we need not make this
assumption. Similar to the discussion at the end of the last
section, we will show that if $\inf_{|x-y|=1} p(x,y)>0$ then the
only possible nonreversible measures are in the case where the
limit of $\pi(x)$ is $0$ in one direction and $\infty$ in the
other direction.

In order to prove the next two propositions we need the following
lemma which appears in a slightly different form as Corollary 5.2
in Liggett(1976):
\begin{lemma}[Liggett] \label{lemma1.1}
If $\inf_{|x-y|=1} p(x,y)>0$ and $\nu\in\mathcal{\tilde{I}}_e$,
then exactly one of the following holds:
\begin{enumerate}
  \item[(a)]
$\nu\{(\eta,\xi):\eta=\xi\}=1$,
  \item[(b)]
$\nu\{(\eta,\xi): \eta\le\xi,\eta\neq\xi\}=1$,
  \item[(c)]
$\nu\{(\eta,\xi):\eta\ge\xi,\eta\neq\xi\}=1$,
  \item[(d)]
$\nu(B)=1$,
  \item[(e)]
$\nu\{(\eta,\xi):(\xi,\eta)\in B\}=1$,
\end{enumerate}
where $B=\{(\eta,\xi):\exists x\in\mathbb{Z}\text{ such that
}\eta(y)\le\xi(y)\text{ for all }y<x, \eta(y)<\xi(y) \text{ for
some }y<x, \eta(z)\ge\xi(z)\text{ for all }z\ge x, \eta(z)>\xi(z)
\text{ for some }z\ge x\}$.
\end{lemma}

\begin{prop}\label{thm1.3}
If $\inf_{|x-y|=1} p(x,y)>0$ and $\pi(x)$ has some finite, nonzero
limit point as $x$ goes to $\infty$ and some finite, nonzero limit
point as $x$ goes to $-\infty$, then $\mathcal{I}_e=\{\nu^{c}:0\le
c\le\infty\}$.
\end{prop}

\begin{proof}
The assumptions imply that $\sum_x{\pi(x)/[1+\pi(x)]^2}=\infty$ so
Theorem \ref{lemma1.2} tells us
$\mathcal{I}_e\supset\{\nu^{c}:0\le c\le\infty\}$. We will show
the reverse containment.

Choose a sequence $\{n_k\}$ extending in both directions so that
finite, nonzero limits of $\pi(n_k)$ exist. For a measure $\mu$ on
$\{0,1\}^{\mathbb{Z}}$ the set of limit points $L_+$ of
$\left\{\mu\{\xi(n_k)=1\}, k>0\right\}$ satisfies one of the
following properties:
\begin{enumerate}
  \item[(i)]
$L_+=\{1\}$ or $L_+=\{0\}$.
  \item[(ii)]
$L_+=\{1,0\}$.
  \item[(iii)] $L_+$ contains some limit point between $0$ and $1$.
\end{enumerate}
The same is true for the set of limit points $L_-$ of
$\left\{\mu\{\xi(n_k)=1\}, k<0\right\}$.

Now suppose we couple $\nu^c$ with another extremal invariant
measure $\mu_e$, the two measures corresponding to the processes
$\eta_t$ and $\xi_t$ respectively.  Since Theorem \ref{lemma1.2}
tells us that $\nu^c$ is extremal, Section VIII.2 in IPS implies
there exists a coupling measure such that $\nu\in
\mathcal{\tilde{I}}_e$.

If $\mu_e$ satisfies condition (i) for both $L_+$ and $L_-$ then
there are two possibilities: either $L_+=L_-$ or $L_+\neq L_-$.
Suppose first that $L_+=L_-=\{1\}$ for $\mu_e$.  If in this case
we have that $\mu_e\{\xi(z)=1\}<1$ for some $z$ then we can choose
$c<\infty$ large enough so that
$\nu^c\{\eta(z)=1\}>\mu_e\{\xi(z)=1\}$. But this contradicts the
assumption that $\nu^c\{\eta(n_k)=1\}=c\pi(n_k)/[1+c\pi(n_k)]$ has
limits less than $1$ for k going to $\infty$ and $-\infty$. To see
this suppose the coupling measure satisfies $\nu(B)=1$ as defined
in Lemma \ref{lemma1.1}.  Given
\begin{equation}
0<\epsilon<1-\lim_{k\rightarrow\infty}c\pi(n_k)/[1+~c\pi(n_k)]
\end{equation}
we can choose $K$ large enough so that
\begin{eqnarray*}
1-\epsilon<\nu\{(\eta,\xi):\exists x<K\text{ such that
}\eta(y)\le\xi(y)\forall y<x, \eta(y)<\xi(y) \text{ for some
}y<x,\\
\eta(z)\ge\xi(z)\forall z\ge x, \eta(z)>\xi(z) \text{ for some
}z\ge x\}.
\end{eqnarray*}
This, however, contradicts $L_+=1$. Similarly we cannot have that
$\nu\{(\eta,\xi):(\xi,\eta)\in B\}=1$. So Lemma \ref{lemma1.1}
tells us that $\eta\le\xi$ which contradicts
$\nu^c\{\eta(z)=1\}>\mu_e\{\xi(z)=1\}$. It must be that
$\mu_e=\nu^\infty$. A similar argument shows that if
$L_+=L_-=\{0\}$ for $\mu_e$ then $\mu_e=\nu^0$.

Consider the second case where $L_-\neq L_+$; without loss of
generality we will assume that $L_-=\{0\}$.

We claim that given $\epsilon>0$, we can find $n$ such that
$\mu_e\{\xi(n)=0\}<\epsilon$ and $\mu_e\{\xi(n+1)=0\}<\epsilon$.
To see this suppose that for some $\epsilon>0$ there exists no $n$
for which this is true. Then since $L_+=\{1\}$, there are
infinitely many $x>0$ for which $\mu_e\{\xi(x)=0\}<\epsilon/4$ and
infinitely many $y>0$ for which $\mu_e\{\xi(y)=0\}\ge\epsilon$.
Choosing $\nu^c$ so that
$\lim_{k\rightarrow\infty}c\pi(n_k)/[1+c\pi(n_k)]=1-\epsilon/2$
gives us a contradiction to Lemma \ref{lemma1.1} and thus proves
the claim.

Given the same $\epsilon>0$ we can choose $m<n$ so that
$\mu_e\{\xi(m-1)=1\}<\epsilon$. Since we have that
$\nu\in\mathcal{\tilde{I}}_e$ then $\int{\tilde{\Omega}(\sum_{x\in
T}{f_x})}d\nu=0$ for each finite ${T}\subset \mathbb{Z}$.  By
(\ref{1}),
\begin{eqnarray}\label{1.5}
&&\sum_{m\le x\le n, y\in\mathbb{Z}}{(p(x,y)+p(y,x))\int f_{yx} d\nu}\\
&=&\sum_{x=m\text{ or }n, y=m-1\text{ or }n+1}{\left[p(x,y)\int
(g_{xy} - h_{yx})d\nu + p(y,x)\int (h_{xy} - g_{yx})
d\nu\right]}\nonumber
\end{eqnarray}
which is increasing in $n$ and $-m$.

Due to our choice of $m$ and $n$ above, $\int
h_{n,n+1}d\nu<\epsilon$ and $\mu_e\{\xi(m-1)=1\}<\epsilon$;
moreover \newline $P(A)-P(A\bigcap B\bigcap C)\le P(B^c)+P(C^c)$
implies that $\nu^c\{\eta(n+1)=1, \eta(n)=0\}-\int
g_{n+1,n}d\nu<2\epsilon$ so that
\begin{eqnarray*}
&&\sum_{m\le x\le n, y\in\mathbb{Z}}{(p(x,y)+p(y,x))\int f_{yx}
d\nu}\\
&<&p(n,n+1)\int g_{n,n+1}d\nu - p(n+1,n)\int g_{n+1,n}d\nu +
3\epsilon\\
&<&p(n,n+1)\nu^c\{\eta(n)=1,
\eta(n+1)=0\}\\
&-&p(n+1,n)\nu^c\{\eta(n)=0, \eta(n+1)=1\}+5\epsilon.
\end{eqnarray*}
By the reversibility of $\nu^c$
\begin{equation*}
p(n,n+1)\nu^c\{\eta(n)=1, \eta(n+1)=0\}=p(n+1,n)\nu^c\{\eta(n)=0,
\eta(n+1)=1\}
\end{equation*}
so equation (\ref{1.5}) is in fact equal to $0$. Since we have
assumed here that $L_-=\{0\}$ and $L_+=\{1\}$ for $\mu_e$, then
choosing $0<c<\infty$ gives us a contradiction.

Suppose $\mu_e$ satisfies condition (ii) for either $L_+$ or $L_-$
so that either $L_+=\{0,1\}$ or $L_-=\{0,1\}$.  Choose $\nu^c$
with $0<c<\infty$. Again we contradict Lemma \ref{lemma1.1}.

Combining all the above arguments we have that either
$\mu_e=\nu^0$, $\mu_e=\nu^\infty$, or $\mu_e$ satisfies (iii) in
some direction.  Assuming the latter we can, without loss of
generality, choose $0<c_0<\infty$ so that
\begin{equation*}
\lim_{k\rightarrow\infty}c_0\pi(n_k)/[1+c_0\pi(n_k)]=\lim_{l\rightarrow\infty}
\mu_e\{\xi(n_{k_l})=1\}.
\end{equation*}
For all $c>c_0$,
\begin{equation*}
\lim_{k\rightarrow\infty}c\pi(n_k)/[1+c\pi(n_k)]>\lim_{l\rightarrow\infty}
\mu_e\{\xi(n_{k_l})=1\}.
\end{equation*}
By Lemma \ref{lemma1.1} either $\mu_e\le \nu^c$ or $\nu(B)=1$
where $B$ is defined in the lemma. Similarly, for all $c<c_0$,
either $\mu_e\ge \nu^c$ or $\nu\{(\eta,\xi):(\xi,\eta)\in B\}=1$.
Combining these two arguments gives
$\nu^{c_1}\le\mu_e\le\nu^{c_2}$ for all $c_1<c_0<c_2$.  By the
continuity of the one parameter family of measures $\nu^c$,
$\mu_e=\nu^{c_0}$.
\end{proof}

\begin{prop}\label{prop1.7}
If $\inf_{|x-y|=1} p(x,y)>0$,
$\lim_{x\rightarrow\infty}\pi(x)=\infty$, and $\pi(x)$ has a
finite, nonzero limit point as $x$ goes to $-\infty$, then
$\mathcal{I}_e=\{\nu^{c}:0\le c\le\infty\}$.
\end{prop}

\begin{proof}
Again, by Theorem \ref{lemma1.2} we need only show that
$\mathcal{I}_e\subset\{\nu^{c}:c\in [0,\infty]\}$.

We argue first that without loss of generality we can assume the
limit points of $\{\pi(x),x<0\}$ are bounded above. Assume to the
contrary that $\infty$ is a limit point.  Then for any $R>0$ we
can find $x<-R$ such that $\min(\pi(x),\pi(x+1))>R$ since
$\inf_{|x-y|=1} p(x,y)>0$. The conditions of Proposition
\ref{thm1.2} are then satisfied so that
$\mathcal{I}_e=\{\nu^{c}:c\in [0,\infty]\}$ holds. We will
therefore assume throughout the rest of the proof that the limit
points of $\{\pi(x),x<0\}$ are bounded above.

Couple $\nu^c$ with another extremal invariant measure $\mu_e$,
the two measures corresponding to the processes $\eta_t$ and
$\xi_t$ respectively. As argued above there exists a coupling
measure such that $\nu\in \mathcal{\tilde{I}}_e$.

Let $L^-$ be the the set of limit points of
$\left\{\mu_e\{\xi(x)=1\}, x<0\right\}$. Note that $L^-$ is
slightly different from $L_-$ described in Proposition
\ref{thm1.3} in that $L_-$ is the set of limit points for a subset
of $\{\mu_e\{\xi(x)=1\}, x<0\}$. $L^-$ satisfies one of the
following properties:
\begin{enumerate}
  \item[(i)]
$L^-$ contains some point between $0$ and $1$.
  \item[(ii)]
$L^-=\{1,0\}$.
  \item[(iii)]
$L^-=\{1\}$.
  \item[(iv)]
$L^-=\{0\}$.
\end{enumerate}
The same is true for the set $L^+$ of limit points of
$\left\{\mu_e\{\xi(x)=1\},x>0\right\}$.

Suppose $L^-$ satisfies (i). Choose a sequence $x_n\rightarrow
-\infty$ so that
$0<\lim_{n\rightarrow\infty}\mu_e\{\xi(x_n)=1\}<1$ exists. Since
we can assume that the limit points of $\{\pi(x),x<0\}$ are all
finite,
 there exists a subsequence $\{x_{n_k}\}$ such that
$\lim_{k\rightarrow\infty} \pi(x_{n_k})<\infty$ exists.

Consider the two cases where $\lim_{k\rightarrow\infty}
\pi(x_{n_k})=0$ and where $\lim_{k\rightarrow\infty}
\pi(x_{n_k})>0$.  Assume the latter case first. Choose
$0<c_0<\infty$ so that
\begin{equation*}
\lim_{k\rightarrow\infty}c_0\pi(x_{n_k})/[1+c_0\pi(x_{n_k})]=
\lim_{n\rightarrow\infty} \mu_e\{\xi(x_{n})=1\}.
\end{equation*}
For all $c>c_0$,
\begin{equation*}
\lim_{k\rightarrow\infty}c\pi(x_{n_k})/[1+c\pi(x_{n_k})]>
\lim_{n\rightarrow\infty} \mu_e\{\xi(x_{n})=1\}.
\end{equation*}
Using the argument at the end of the proof of Proposition
\ref{thm1.3}, we have that for all $c_1<c_0<c_2$,
$\nu^{c_1}\le\mu_e\le\nu^{c_2}$. Consequently, it must be that
$\mu_e=\nu^{c_0}$.

Now assume that $\lim_{k\rightarrow\infty} \pi(x_{n_k})=0$ so that
for all $0<c<\infty$ the coupling satisfies either $\nu^c\le\mu_e$
or $\nu\{B\}=1$ where $B$ is given in Lemma \ref{lemma1.1}. If
$\nu^c\le\mu_e$ for all $0<c<\infty$ then $\mu_e=\nu^\infty$, a
contradiction to $L^-$ satisfying (i). So it must be that
$\nu\{B\}=1$.

We claim that for any $r<1$ there exists $m<0$ such that
$\mu_e\{\xi(m)=1\}>r$ and $\mu_e\{\xi(m-1)=1\}>r$.  By the
hypothesis of the theorem we can choose a sequence $\{x_l\}$ going
to $-\infty$ so that $
 0<\lim_{l\rightarrow\infty}
\pi(x_l)<\infty$ exists.  If $\inf_{|x-y|=1} p(x,y)>p$
 then choose
$c$ so that
\begin{equation*}
\lim_{l\rightarrow\infty}\frac{cp\pi(x_l)}{1+cp\pi(x_l)}>r+\frac{1-r}{2}.
\end{equation*}
Since $\pi(x_l-1)>p\pi(x_l)$, it follows that the set of limit
points of $\{\frac{c\pi(x_l-1)}{1+c\pi(x_l-1)}, l>0\}$ is bounded
below by $r+\frac{1-r}{2}$. Now since $\nu\{B\}=1$ there exists a
$K$ such that $l>K$ implies $\mu_e\{\xi(x_l)=1\}>r$ and
$\mu_e\{\xi(x_l-1)=1\}>r$ which proves the claim.

Since we have that $\nu\in\mathcal{\tilde{I}}_e$ then
$\int{\tilde{\Omega}(\sum_{x\in T}{f_x})}d\nu=0$ for each finite
${T}\subset \mathbb{Z}$.  By (\ref{1}),
\begin{eqnarray*}
&&\sum_{m\le x\le n, y\in\mathbb{Z}}{(p(x,y)+p(y,x))\int f_{yx} d\nu}\\
&=&\sum_{x=m\text{ or }n, y=m-1\text{ or }n+1}{\left[p(x,y)\int
(g_{xy} - h_{yx})d\nu + p(y,x)\int (h_{xy} - g_{yx})
d\nu\right]}\nonumber
\end{eqnarray*}
which is increasing in $n$ and $-m$.

Using the claim above along with the fact that
$\lim_{x\rightarrow\infty}\pi(x)=\infty$, we can argue just as we
argued in the case where $L_-\neq L_+$ of (i) in Proposition
\ref{thm1.3}, to get
\begin{eqnarray*}
&&\sum_{m\le x\le n, y\in\mathbb{Z}}{(p(x,y)+p(y,x))\int f_{yx}
d\nu}\\
&<&p(m,m-1)\int g_{m,m-1}d\nu - p(m-1,m)\int g_{m-1,m}d\nu +
3\epsilon\\
&<&p(m,m-1)\nu^c\{\eta(m)=1,
\eta(m-1)=0\}\\
&-&p(m-1,m)\nu^c\{\eta(m)=0, \eta(m-1)=1\}+5\epsilon.
\end{eqnarray*}
By the reversibility of $\nu^c$ the left-hand side is just
$5\epsilon$, but this contradicts $\nu\{B\}=1$ for small
$\epsilon$.

Suppose $L^-$ satisfies condition (ii). Choosing $\nu^c$ with
$0<c<\infty$ gives us a contradiction to Lemma \ref{lemma1.1}.

If $L^-$ satisfies condition (iii) then we will handle the two
cases (a) $L^+=\{1\}$ and (b) $L^+\neq \{1\}$. Consider case (a)
first. If we switch the coupling so that $\mu_e$ corresponds to
$\eta_t$ and $\nu^c$ corresponds to $\xi_t$ then we have that the
left-hand side of the following inequality goes to $0$:
\begin{eqnarray}\label{eqn9}
&&\sum_{|x|=n, |y|=n+1}{(p(x,y)+p(y,x))\int f_yd\nu}\ge\\
&&\sum_{|x|=n, |y|=n+1}{p(x,y)\int (g_{xy} - h_{yx})d\nu}
+\sum_{|x|=n, |y|=n+1}{p(y,x)\int (h_{xy} - g_{yx}) d\nu}\nonumber
\end{eqnarray}
By (\ref{1}) and by irreducibility we get $\int f_{xy} d\nu=0$ for
all $x,y$. Therefore the measure $\mu_e$ must lie stochastically
above all $\nu^c$ for all finite $c$ and must therefore be equal
to $\nu^{\infty}$.

If (b) holds then we refer the reader to the argument given above
in the case where $L^-$ satisfies (i) and
$\lim_{k\rightarrow\infty} \pi(x_{n_k})=0$.

Finally suppose that (iv) holds so that $L^-=\{0\}$.  If $L^+$
satisfies (i) or (ii) then by Lemma \ref{lemma1.1}, $\mu_e\le
\nu^c$ for all $c>0$ so that $\mu_e=\nu^0$, a contradiction. If
$L^+$ satisfies (iv) then similarly $\mu_e=\nu^0$. Let $L^+$
satisfy (iii) so that $L^+=\{1\}$. For a given $z$ choose $c$
small enough so that $\nu^c\{\eta(z)=1\}<\mu_e\{\xi(z)=1\}$. We
thus have that $\nu\{(\eta,\xi):(\xi,\eta)\in B\}=1$ as given in
Lemma \ref{lemma1.1}.  But by (\ref{1}) and (\ref{eqn9}), for a
given $\epsilon>0$ we can find $-m$ and $n$ large enough so that
\begin{equation*}
\sum_{m\le x\le n, y\in\mathbb{Z}}{(p(x,y)+p(y,x))\int
f_{yx}d\nu}<\epsilon
\end{equation*} which of course contradicts
$\nu\{(\eta,\xi):(\xi,\eta)\in B\}=1$.
\end{proof}

\begin{proof}[Proof of Theorem \ref{thm1.8}]
Note first that since $\inf_{|x-y|=1} p(x,y)>0$ then it cannot be
that $\mathcal{L}^-$ or $\mathcal{L}^+$ is equal to
$\{0,\infty\}$. In light of this fact, if either $\mathcal{L}^-$
or $\mathcal{L}^+$ contains a finite, nonzero point then
Proposition \ref{thm1.3} and analogs of Proposition \ref{prop1.7}
imply there are no nonreversible measures.  If
$\mathcal{L}^+=\mathcal{L}^-=\{0\}$ or
$\mathcal{L}^+=\mathcal{L}^-=\{\infty\}$ then Proposition
\ref{thm1.2} implies there are no nonreversible measures.
\end{proof}

\section{A result concerning domains of attraction}\label{domains}
\begin{theorem}\label{thmlastsection}
Let $\sum_x \pi(x)/[1+\pi(x)]^2=\infty$ and let $\omega$ be a
probability measure on $[0,\infty]$. Also, assume that $\nu_c$ is
a family of invariant measures indexed by $c\ge 0$ each of which
is in $\mathcal{I}_e$. Suppose $\{\mu_c\}$ is a family of
probability measures on $\{0,1\}^{\mathcal{S}}$ such that for each
$0\le c\le\infty$, $\mu_c$ is absolutely continuous with respect
to $\nu_c$. If
\begin{equation}\label{defmu}
\mu=\int_0^\infty \mu_c \,\omega(dc)\text{ and }\nu=\int_0^\infty
\nu_c \,\omega(dc)
\end{equation}
then $ \lim_{T\rightarrow\infty}\frac{1}{T}\int_0^T\mu S(t) dt$
exists and is equal to $\nu$.
\end{theorem}

\begin{proof}
 For a fixed $c$ we first prove that
\begin{equation}\label{equ1}
\lim_{T\rightarrow\infty}\frac{1}{T}\int_0^{T} \mu_c S(t)
dt=\nu_c.
\end{equation}

By the compactness of $\mathcal{P}$ we can choose a sequence of
times such that
\begin{equation}\label{equ2}
\lim_{n\rightarrow\infty}\frac{1}{t_n}\int_0^{t_n} \mu_cS(t) dt
\end{equation} converges in distribution
to some measure $\lambda$. Pick a continuous (and therefore
bounded) function $f$ on $\{0,1\}^\mathcal{S}$ with $\|f\|\le 1$
and let $g$ be the Radon-Nikodym derivative of $\mu_c$ with
respect to $\nu_c$. Given $\epsilon>0$ we have that for $n$ large
enough
\begin{equation*}
|\frac{1}{t_n}\int_0^{t_n} \int (S(t) f)g\,d\nu_c dt-\int
f\,d\lambda|<\epsilon/3.
\end{equation*}

We can choose a simple function
\begin{equation*}
\hat{g}=\sum_{k=1}^N c_k1_{E_k}
\end{equation*}
 approximating $g$ such that
$\cup_k E_k=\{0,1\}^\mathcal{S}$, $\hat{g}\ge 0$, $\int \hat{g}\,
d\nu_c=1$, and $\int |g-\hat{g}|\,d\nu_c<\epsilon/3$. Since
$\|S(t)f\|\le\|f\|\le 1$ this gives us
\begin{equation*}
|\frac{1}{t_n}\int_0^{t_n} \int (S(t) f)g\,d\nu_c dt-
\frac{1}{t_n}\int_0^{t_n} \int (S(t) f)\hat{g}\,d\nu_c dt|\le \int
|g-\hat{g}|\,d\nu_c<\epsilon/3.
\end{equation*}

Without loss of generality we can henceforth assume that
$\nu_c(E_k)>0$ for each $k$. Define the measure $\mu_k$
concentrating on $E_k$ by letting
\begin{equation*}
\mu_k(A)=\frac{\nu_c(A)}{\nu_c(E_k)}
\end{equation*}
 for all $A\subset E_k$ and
$\mu_k=0$ otherwise. If we think of $\hat{g}$ as the Radon-Nikodym
derivative of some measure $\lambda_{\epsilon}$ with respect to
$\nu_c$ then we can write
\begin{equation*}
\sum_{k=1}^N \nu_c(E_k)\mu_k=\nu_c\text{ and }\sum_{k=1}^N
c_k\nu_c(E_k)\mu_k=\lambda_{\epsilon}.
\end{equation*}

We can now find a subsequence $\{t_{n_l}\}$ such that the
following limits exist for each $k$:
\begin{equation*}
\lim_{l\rightarrow\infty}\frac{1}{t_{n_l}}\int_0^{t_{n_l}}
\mu_kS(t) dt=\nu_k.
\end{equation*}
 Moreover, Proposition I.1.8 in IPS tells us
$\nu_k\in\mathcal{I}$. Since $\nu_c$ is extremal invariant and
since $\sum_{k\ge 1} \nu_c(E_k)\nu_k=\nu_c$, it must be that
$\nu_k=\nu_c$ for each $k$. This then yields
\begin{equation*}
\sum_{k=1}^N
c_k\nu_c(E_k)\nu_k=\lim_{l\rightarrow\infty}\frac{1}{t_{n_l}}\int_0^{t_{n_l}}
\lambda_{\epsilon}S(t) dt=\nu_c
\end{equation*}
which gives us
\begin{equation*}
|\frac{1}{t_{n_l}}\int_0^{t_{n_l}} \int (S(t) f)\hat{g}\,d\nu_c
dt-\int f\,d\nu_c|<\epsilon/3
\end{equation*}
for $l$ large enough.

Combining the three inequalities we have
\begin{equation*}
|\int f\,d\lambda-\int f\,d\nu_c|<\epsilon.
\end{equation*}
But $\epsilon>0$ is arbitrary so it must be that $\int
f\,d\lambda=\int f\,d\nu_c$ for each continuous $f$ with $\|f\|\le
1$ which implies that (\ref{equ2}) is equal to $\nu_c$. Now let
$M_n$ be the closure of the set of measures
\begin{equation*}\{\frac{1}{T}\int_0^{T} \mu_cS(t) dt:T\ge n\}.
\end{equation*}
Using the compactness of $\mathcal{P}$ along with the fact that
$\{t_n\}$ is an arbitrary sequence of times causing convergence in
(\ref{equ2}), we have that $\bigcap_{n\in\mathbb{N}} M_n=\nu_c$
proving (\ref{equ1}).

To finish the proof note that since $\|S(t)f\|\le\|f\|$, we can
use the Dominated Convergence Theorem together with Fubini's
Theorem to show that
\begin{equation*}
\lim_{T\rightarrow\infty}\frac{1}{T}\int_0^{T} \int_0^\infty \int
S(t)f \,d\mu_c \,\omega(dc)
 \,dt=\int f\, d\nu.
\end{equation*}
\end{proof}

For the following corollary let $\nu_\alpha$ be the product
measure with marginals
$0<\nu_\alpha\{\eta:\eta(x)=1\}=\alpha(x)<1$ for $\alpha(x)$ a
function on $\mathcal{S}$.
\begin{corollary}
Suppose $\sum_x \pi(x)/[1+\pi(x)]^2=\infty$.  If $\sum_x
|\alpha(x)-\frac{c\pi(x)}{1+c\pi(x)}|<\infty$ then
\begin{equation}\label{convergence}
\lim_{T\rightarrow\infty}\frac{1}{T}\int_0^T\nu_\alpha S(t)
dt=\nu^c.
\end{equation}
\end{corollary}
\begin{proof}
Let $\beta(x)=\frac{c\pi(x)}{1+c\pi(x)}$,
$m_x=\min[\alpha(x),\beta(x)]$, and
$M_x=\max[\alpha(x),\beta(x)]$. We then have
\begin{eqnarray*}
1-|\alpha(x)-\beta(x)|&=&1-M_x+m_x \\
&=&
[(1-M_x)(1-M_x)]^{1/2}+(m_xm_x)^{1/2}\\
&\le&
[(1-M_x)(1-m_x)]^{1/2}+(m_xM_x)^{1/2}\\
&=&[(1-\alpha(x))(1-\beta(x))]^{1/2}+(\alpha(x)\beta(x))^{1/2}.
\end{eqnarray*}
Since $\sum_x |\alpha(x)-\beta(x)|<\infty$ then
\begin{equation*}
\prod_x
\{(\alpha(x)\beta(x))^{1/2}+[(1-\alpha(x))(1-\beta(x))]^{1/2}\}\ge
\prod_x \{1-|\alpha(x)-\beta(x)|\}>0.
\end{equation*}
An application of Kakutani's Dichotomy tells us that $\nu_\alpha$
is absolutely continuous with respect to $\nu^c$ which completes
the proof.
\end{proof}
We remark here that if $\alpha(x)$ and $\beta(x)$ are both bounded
away from $0$ and $1$ then Kakutani's Dichotomy tells us that
$\sum_x [\alpha(x)-\beta(x)]^2<\infty$ is a necessary and
sufficient condition for $\nu_\alpha$ to be absolutely continuous
with respect to $\nu^c$ (e.g. page 245 of Durrett(1996)).

\textbf{Acknowledgement}. The author thanks his advisor, Thomas M.
Liggett, for motivating this paper and for the many discussions
that led to the writing of it.

\end{document}